\documentclass[11pt]{amsart}
\usepackage{latexsym}
\usepackage{amssymb,amscd,amsthm,amsfonts,verbatim,amsopn,color}
\usepackage[dvips]{epsfig}
\def\Z{{\mathbb Z}}

\def\C{{\mathbb C}}
\def\F{{\mathbb F}}
\def\R{{\mathbb R}}

\def\P{{\mathbb P}}

\def\AA{{\mathcal A}}

\def\CC{{\mathcal C}}
\def\DD{{\mathcal D}}
\def\GG{{\mathcal G}}
\def\LL{{\mathcal L}}
\def\MM{{\mathcal M}}
\def\NN{{\mathcal N}}

\def\Ss{{\mathcal S}}
\def\TT{{\mathcal T}}

\newcommand{\aut}{{\rm Aut}}
\newcommand{\lra}{\longrightarrow}
\newcommand{\stab}{{\rm stab}}
\newcommand{\mfB}{\mathfrak B}
\newcommand{\mfS}{\mathfrak S}

\newcommand{\mfD}{\mathfrak D}
\newcommand{\ol}{\overline}
\newcommand{\proj}{\text{proj}}
\newcommand{\ra}{\rightarrow}
\newcommand{\ti}{\tilde}
\newcommand{\vrt}{\text{vert}\,}
\newcommand{\wti}{\widetilde}
\newcommand{\YP}{Y_{\Pi_n}}

\newtheorem{thm}{Theorem}[section]
\newtheorem{df} [thm]{Definition}
\newtheorem{lm}  [thm]{Lemma}

\newtheorem{prop}[thm]{Proposition}

\newtheorem{rem}[thm]{Remark}

\newtheorem{expl}[thm]{Example}
\numberwithin{equation}{section}

\newenvironment{pf}{\noindent {\bf Proof.}}{\hfill $\Box$\vspace{0.3cm}}

\begin{document}

\title[Abelianizing the real permutation action via blowups]
{Abelianizing the real permutation action \\ via blowups}

\author{Eva Maria Feichtner \mbox{ }\& \mbox{ }Dmitry N.\ Kozlov}

\address{
Department of Mathematics, ETH Zurich, 8092 Zurich, Switzerland}
\email{feichtne@math.ethz.ch}

\address{
Department of Mathematics, University of Bern, 3012 Bern, Switzerland;
\newline
\mbox{ }\hspace{0.1cm}\,
{\em on leave from}: Department of Mathematics, KTH Stockholm,
100 44 Stockholm, Sweden}
\email{kozlov@math.kth.se}

\date{February 2003}

%
%
%

\maketitle

\section{Introduction}

Our object of study is an abelianization of the $\Ss_n$ permutation action
on $\R^n$ that is provided by a particular De~Concini-Procesi wonderful 
model for the braid arrangement.
Our motivation comes from an analogous construction for
finite group actions on complex manifolds, due to Batyrev \cite{B1,
  B2}, and subsequent study of Borisov \& Gunnells \cite{BG}, where
the connection of such abelianizations with De~Concini-Procesi
wonderful models for arrangement complements was first observed.

Whereas previous studies were restricted to complex manifolds, here we
study one of the most natural nontrivial actions of a~finite group on
a~{\em real\/} differentiable manifold, namely the permutation action
on~$\R^n$. The locus of non-trivial stabilizers in this case is
provided by the braid arrangement~$\AA_{n-1}$. We suggest to blow up
intersections of subspaces in~$\AA_{n-1}$, respectively proper
transforms of those intersections, in the order of an~arbitrary linear
extension of the intersection lattice~$\Pi_n$, so as to exhaust all of
the arrangement.  That is the same as to take the De~Concini-Procesi
wonderful model of the arrangement complement with respect to the 
maximal building set, see \cite{DP}.

Not only do we obtain an abelianization of the real permutation
action, we even show that stabilizers of points in the arrangement
model are isomorphic to direct products of~$\Z_2$. To this end, we
develop a~combinatorial framework for explicitly describing the
stabilizers in terms of automorphism groups of set diagrams over
families of cubes.

Moreover, we observe that the natural nested set stratification on the
arrangement model is not stabilizer distinguishing with respect to the
$\Ss_n$-action, i.e., stabilizers of points are not in general
isomorphic on open strata.  Motivated by this structural deficiency,
we furnish a~new stratification of the De~Concini-Procesi arrangement
model that distinguishes stabilizers.

Arrangement models have been extensively studied over the last years.
They were introduced by De~Concini~\& Procesi in \cite{DP}, one of the
motivations being to provide rational models for cohomology algebras
of arrangement complements. In \cite{FK} the De~Concini-Procesi model
construction was put in a~very general combinatorial context, showing
that the notions of building sets and nested sets, coined already by
Fulton~\& MacPherson in~\cite{FM}, along with the notion of a~blowup,
have canonical combinatorial counterparts in the theory of
semilattices. It was also shown in \cite{FK} that this combinatorial
framework actually traces precisely the step-by-step change in the
incidence structure of strata during the De~Concini-Procesi resolution
process.

On the geometric side, wonderful arrangement models were generalized
to wonderful conical compactifications by MacPherson~\&
Procesi~\cite{MP}, and Gaiffi~\cite{G2} recently provided a further
generalization incorporating mixed real subspace and halfspace
arrangements as well as real stratified manifolds as starting points
of the construction. Algebraic topological invariants of wonderful
models are another focus of interest. Yuzvinsky~\cite{Y} provided a
monomial basis for the cohomology of wonderful compactifications of
hyperplane arrangements that was later generalized by Gaiffi to
compactifications of subspace arrangements in \cite{G1}.

We give a more detailed outline of our paper: In Section~\ref{sect_DPmodels}, 
we begin our investigations with a brief review of De~Concini-Procesi
wonderful models. Moreover, we describe
how an action of a finite group on an arrangement extends to an action
on the arrangement model.  We then turn to our specific situation,
observing that when blowing up the entire locus of non-trivial
stabilizers for $\Ss_n$ acting on $\R^n$, i.e., the entire braid
arrangement, the nested set stratification is not sufficient to
distinguish stabilizers. That is, we may have two points lying on
the same stratum, but having non-isomorphic stabilizers. In fact, this
happens already for $n=3$.

In Section~\ref{sect_nstrat}, we study the nested set
stratification and group actions on De~Concini-Procesi models in
some detail, so that finally, in Section \ref{sect_stdist}, we are
able to rectify the situation: We define a~different
stratification on the De~Concini-Procesi model such that, on one
hand, this stratification is naturally arrived at by tracing
a~certain, interesting on its own right, subspace arrangement in
$\R^n$, on the other hand, this new stratification is {\em
stabilizer distinguishing\/}.

In Section \ref{sect_combfr} we turn to the detailed study of the
isomorphism types of stabilizers of points in the De~Concini-Procesi
resolution of the braid arrangement. Relying on our analysis in the
previous sections, we know that the stabilizer of a~point in the
arrangement model is the intersection of a number of stabilizers of
lines and of the stabilizer of one single point in $\R^n$.  We develop
a combinatorial language to describe stabilizers of points and lines
in $\R^n$, namely by representing them as automorphism groups of set
diagrams over families of cubes. The crucial property of this
representation is that taking intersections of a number of
automorphism groups of such diagrams will again yield an~automorphism
group over a diagram. This new diagram can be combinatorially read of
from the original diagrams. Thus, we succeed to represent the
stabilizer of a point in the arrangement model as an automorphism
group of a set diagram over a family of cubes. By further analysis of
this diagram, we are finally able to prove in
Section~\ref{sect_stabYP} that, beyond the natural initial expectation
that the stabilizers ought to be abelian, they in fact are isomorphic
to direct products of~$\Z_2$, with the number of factors in
each product at most $\lfloor\frac{n}{2}\rfloor$.


\section{De Concini-Procesi arrangement models}
\label{sect_DPmodels}

In this section we briefly review the construction and main
characteristics of wonderful arrangement models as introduced by
De~Concini \& Procesi in~\cite{DP}. We first remind the notions of
building sets and nested sets since they guide the explicit
construction and capture the underlying incidence combinatorics of
a~natural stratification. Moreover, we comment on actions of finite
groups on De~Concini-Procesi models that are induced from group actions on
the arrangement.

\subsection{Building sets and nested sets}

Let $\AA$ be an arrangement of linear subspaces in a~finite
dimensional real or complex vector space, and denote by
$\LL=\LL(\AA)$ the lattice of intersections of spaces in $\AA$
ordered by reverse inclusion, customarily called the {\bf
intersection lattice} of $\AA$.

\begin{df} \label{df_DPnotions}
{\rm (\cite[\S 2]{DP})} For $\LL\,{=}\,\LL(\AA)$ the intersection
lattice of a complex or real subspace arrangement, let $\LL^*$
denote the lattice formed by the orthogonal complements of
intersections in $\AA$ ordered by inclusion.
\begin{itemize}
\item[(1)] For $U\,{\in}\,\LL^*$,
$U\,{=}\,\oplus_{i=1}^k U_i$ with $U_i\,{\in}\,\LL^*$, is called a
{\bf decomposition\/} of~$U$ if for any $V\,{\subseteq}\,U$,
$V\,{\in}\, \LL^*$, $V\,{=}\,\oplus_{i=1}^k (U_i\,{\cap}\,V)$ and
$U_i\,{\cap}\, V\,{\in}\,\LL^*$, for $i=1,\ldots,k$.
\item[(2)] Call  $U\,{\in}\,\LL^*$ {\bf irreducible\/} if it does not
admit a non-trivial decomposition.
\item[(3)] $\GG\,{\subseteq}\,\LL^*\,{\setminus}\,\{\hat 0\}$ is called a
{\bf building set\/}
for~$\AA$ if for any $U\,{\in}\,\LL^*\,{\setminus}\,\{\hat 0\}$ 
and $G_1,\ldots, G_k$ maximal in~$\GG$ below~$U$,
$U\,{=}\, \oplus_{i=1}^k G_i$ is a decomposition (the $\GG$-decomposition)
of~$U$.
\item[(4)] A subset $\TT\,{\subseteq}\,\GG$ is called {\bf nested\/} if
for any set of non-comparable elements $U_1,\ldots,U_k$ in~$\TT$,
$U\,{=}\,\oplus_{i=1}^k\, U_i$ is the $\GG$-decomposition of~$U$. The
nested sets in $\GG$ form an abstract simplicial complex, the
{\bf nested set complex} $\NN(\GG)$.
\end{itemize}
\end{df}

We will without further notice consider building sets as subsets of
the intersection lattice $\LL$, and thus let the consideration of $\LL^*$
remain a detour for the sake of providing a transparent definition. Note that
for any arrangement $\AA$ the set of irreducible elements in 
$\LL(\AA)\,{\setminus}\,\{\hat 0\}$ is
the minimal building set, whereas 
$\GG\,{=}\,\LL(\AA)\,{\setminus}\,\{\hat 0\}$     is the
maximal building set. For the maximal building set the nested
set complex coincides with the order complex of the (non-reduced)
intersection lattice.

\subsection{Arrangement models and the nested set stratification}

We are now prepared to give the definition of wonderful
arrangement models. Let~$\AA$ be an arrangement of subspaces in a
real or complex vector space~$V$, $\LL(\AA)$ its intersection
lattice, and $\GG$ a building set for~$\AA$. On the complement of
the arrangement, $\MM(\AA)\,{:=}\, V\,{\setminus}\,\bigcup \AA$,
consider the map
\begin{equation}\label{eqn_defmapYG}
    \Phi: \quad \MM(\AA) \, \, \longrightarrow \, \,
                    V\, \times \, \prod_{G\in\GG}\, \P(V/G)\, ,
\end{equation}
where in its first coordinate the map is given by inclusion, and in
later coordinates by projection to the (real, resp.\ complex)
projectivizations of the respective quotient spaces. Formally,
\[
\Phi(x)\, \, =\, \, (\,x\,,\,(\Phi_G(x))_{G\in\GG}\,)\,,
\]
with  $\Phi_G(x)=\langle x,G\rangle/G\in\P(V/G)$, for $x\in
\MM(\AA)$, where brackets $ \langle \cdot\, , \cdot \rangle$ denote 
the linear span of subspaces or vectors, respectively. 
This map is an embedding of $\MM(\AA)$, the arrangement
model $Y_{\GG}$ is defined as the closure of its image in $V\,
\times \, \prod_{G\in\GG}\, \P(V/G)$:
\[
               Y_{\GG} \, \, :=\, \, {\rm cl}\, ({\rm Im}\, \Phi)\, .
\]

Alternatively, $Y_{\GG}$ can be described as the result of
subsequently blowing up intersections of subspaces in $\AA$, and
proper transforms of such, corresponding to building set elements
$G\in \GG$ in some linear extension of the inclusion order.

The arrangement model~$Y_{\GG}$ is a smooth variety that contains
the arrangement complement~$\MM(\AA)$ as an open subspace. The
complement $D$ of $\MM(\AA)$ in~$Y_{\GG}$ is a divisor with normal
crossings, in fact, it is the union of smooth, irreducible
components $D_G$ indexed by building set elements~$G\,{\in}\,\GG$.
The intersections of divisors $D_G$ are smooth and irreducible,
naturally, they are indexed with subsets of $\GG$. One of the main
results of De Concini and Procesi, \cite{DP}, states that an
intersection of divisors is non-empty if and only if it is indexed
with a~nested set in~$\GG$.

We call the resulting stratification of $Y_{\GG}$ by irreducible
divisor components $D_{G}$ and their intersections the {\em nested
set stratification\/} of $Y_{\GG}$, and denote it by
$(Y_{\GG},\mfD)$. Note that the poset of strata for
$(Y_{\GG},\mfD)$ coincides with the face poset of the nested set
complex~$\NN(\GG)$.

De~Concini \& Procesi also provide a projective version of their
arrangement models obtained by starting out with the projectivization
of the arrangement complement and replacing the first factor on the
right hand side of~(\ref{eqn_defmapYG}) by~$\P(V)$ accordingly. The
properties of the resulting projective model $\ol Y_{\GG}$ are
similar to those of~$Y_{\GG}$, for details we refer to~\cite[\S
4]{DP}.

\subsection{Finite group actions on arrangements and on their 
wonderful models}
\label{ssect_gract}

Let us now assume that a~finite group $\Gamma$ acts on our vector
space $V$ by linear transformations, and that the arrangement $\AA$ is
invariant under that action. By a standard result from representation
theory, any linear action of a finite group is orthogonal~\cite[2.3,
Thm.~1]{V}. Throughout the paper, we denote the corresponding 
$\Gamma$-invariant positive definite symmetric bilinear form by the 
usual scalar product.

Since we assume~$\Gamma$ to preserve~$\AA$, the group acts on the 
intersection lattice of $\AA$, 
\[
\gamma(A_1\cap\dots\cap A_r)\, \, = \, \, 
\gamma(A_1)\cap\ldots\cap\gamma(A_r)\,, \quad \mbox{ for all}\,
\gamma\,{\in}\,\Gamma,\, \, A_1,\dots,A_r\,{\in}\,\AA\, , 
\]
as well as internally on the corresponding intersections of subspaces.
Also, $\Gamma$ acts on the ambient space of the arrangement model
corresponding to the maximal building set, that is on
$V\,{\times}\,\prod_{G\in\GG}\P(V/G)$, where 
$\GG=\LL(\AA)\,{\setminus}\,\{\hat 0\}$, by
\begin{eqnarray*}
   \gamma \,(x,(x_G)_{G\in\GG}) & = &
   (\,\gamma(x), (\,\gamma(x_{\gamma^{-1}(G)})\,)_{G\in\GG}), \\
& & \qquad \mbox{ for all }\,\gamma\in\Gamma,\,\,\,
   (\,x,(x_G)_{G\in \GG}) \in V\times \prod_{G\in\GG}\,\P(V/G)\, .
\end{eqnarray*}

Moreover, the inclusion map $\Phi:\,\MM(\AA)\,{\longrightarrow}\,
V\,{\times}\,\prod_{G\in\GG}\P(V/G)$ defined in~(\ref{eqn_defmapYG})
commutes with the action of $\Gamma$:
\begin{eqnarray*}
\gamma(\Phi(x)) & = & \gamma \,(\,x, (\,\langle
x,G\rangle/G\,)_{G\in \GG}) \, \, = \, \, (\,\gamma(x),
     (\,\gamma(\,\langle x,\gamma^{-1}(G)\rangle/\gamma^{-1}(G))\, )
           _{G\in \GG})  \\
                & = &
      (\,\gamma(x),
        (\,\langle \gamma(x),G\rangle/G)_{G\in \GG} )  \, \, = \, \,
                    \Phi(\gamma(x)), \quad \mbox{ for }\,
                            \gamma \in \Gamma, \,
                             x\in \MM(\AA)\, .
\end{eqnarray*}
We conclude that, since each element of $\Gamma$ acts continuously
on $V$, the closure of Im$\,\Phi$ is $\Gamma$-invariant. Hence,
$\Gamma$ acts on the arrangement model $Y_{\GG}$ extending the 
$\Gamma$-action on $\MM(\AA)\,{\subseteq}\,Y_{\GG}$.

Note that choosing a $\Gamma$-invariant building set 
$\GG\,{\subsetneq}\,\LL(\AA)\,{\setminus}\,\{\hat 0\}$ 
as well yields an action of~$\Gamma$ on the corresponding arrangement model.


\section{The arrangement model $Y_{\Pi_n}$}

\subsection{A candidate for an abelianization of the permutation action}

We consider the permutation action of the symmetric group~$\Ss_n$ on $\R^n$,
\[
\qquad \sigma(x)\, \, = \, \,
                      (x_{\sigma(1)},\ldots,x_{\sigma(n)}), \qquad
             \mbox{ for all }\, \sigma \in \Ss_n, \, \,
                                x=(x_1,\ldots,x_n)\in \R^n\, .
\]
The locus of points in $\R^n$ with non-trivial stabilizer is a
union of hyperplanes~$H_{i,j}$, $H_{i,j}\,{:=}\,$ker$(x_i{-}x_j)$
for $1{\leq} i{<}j{\leq}n$. This family of ``diagonal
hyperplanes'' in $\R^n$ is the {\em
  braid arrangement} $\AA_{n{-}1}$ of rank~$n{-}1$, its name referring
to the fact that the complement of a~complexified version in
$\C^n$ is the classifying space of the pure braid group on $n$
strands. The braid arrangement is one of the central examples in
arrangement theory and has provided a starting point for many
investigations and developments in arrangement theory and beyond,
see e.g.,~\cite{OT}.

The intersection lattice of $\AA_{n{-}1}$ is the {\em partition
lattice}~$\Pi_n$, i.e., the poset of set partitions
$\pi=(\pi_1|\ldots|\pi_r)$ of $\{1,\ldots, n\}=:[n]$,
$\pi_i\,{\subseteq}\,[n]$ with $\bigcup_{i=1}^r\pi_i=[n]$, ordered
by reverse refinement. Clearly, a partition
$\pi=(\pi_1|\ldots|\pi_r)$ in $\Pi_n$ corresponds to the
intersection of hyperplanes $\bigcap_{(i,j)\in J_{\pi}}H_{i,j}$
with $J_{\pi}\,{=}\,\{(i,j)\,|\,1{\leq}i{<}j{\leq}n,\,
\{i,j\}\,{\subseteq}\,\pi_k$, for some $1{\leq}k{\leq}r\}$. We
will freely use this correspondence between partitions and
intersections of subspaces in the braid arrangement.

For further considerations, we restrict the permutation action to
the $(n{-}1)$-dimensional real space
\[
    V \, \, = \, \, \{x\in \R^n\,|\, \sum_{i=1}^n\, x_i=0\,\}\, .
\]
The locus of points in $V$ with non-trivial stabilizers is the
intersection of $\AA_{n{-}1}$ with $V$, an essential arrangement with
intersection lattice~$\Pi_n$, which we still call braid arrangement and
denote by $\AA_{n{-}1}$ without further mention.

We propose to study the De~Concini-Procesi arrangement model
$Y_{\Pi_n}$ for $\AA_{n{-}1}$ as a candidate for an abelianization
of the permutation action. 
We allow ourselves here to use the shorthand notation $Y_{\Pi_n}$
instead of $Y_{\Pi_n\,{\setminus}\,\{\hat 0\}}$.
It follows from the general discussion
in subsection \ref{ssect_gract} that $Y_{\Pi_n}$
carries a natural $\Ss_n$-action extending the $\Ss_n$-action on
$\MM(\AA_{n{-}1})\subseteq Y_{\Pi_n}$. It turns out that rather
curious phenomena enter the scene already in low dimensions.

\subsection{The nested set stratification is not stabilizer
distinguishing}
\label{ssect_chebpnts}

Already for $S_3$ acting on $\R^3$, the nested set stratification
on the De~Concini-Procesi model, $(Y_{\Pi_3},\mfD)$, is not fine
enough to distinguish stabilizers. Let us have a close look at the
situation.

As above, we restrict the permutation action to
$V\,{=}\,\{\,(x_1,x_2,x_3)\,|\,
                      \sum_{i=1}^3x_i\,{=}\,0\,\}\,{\subseteq}\,\R^3$.
The arrangement model $Y_{\Pi_3}$ is the result of blowing up
$\{0\}$ in~$V$. Topologically, $Y_{\Pi_3}$ is an open M\"obius
band. As a subspace of $V\,{\times}\,\P(V)$, $Y_{\Pi_3}$ can be
described as follows:
\[
Y_{\Pi_3}\, = \, \{\,(x,\langle x \rangle)\,|\, x\neq 0\,\} \, \cup \,
                 \{\,(0,l)\,|\, l\in \P(V)\,\}
         \, \subseteq \, V \times \P(V)\, .
\]
In terms of this pointwise description of $Y_{\Pi_3}$ the divisors
$D_G$, $G\,{\in}\,\Pi_3$, read
\begin{eqnarray*}
D_{\{0\}}\,=\,D_{(1,2,3)}
             & = &  \{\,(0,l)\,|\, l\in \P(V)\,\} \\
D_{(1,2)(3)} & = & \{\,(x,\langle x \rangle)\,|\, x_1=x_2\neq 0\,\} \,
\cup \, \{\,(0,\langle (1,1,-2)\rangle)\,\} \, ,
\end{eqnarray*}
with $D_{(1,3)(2)}$, $D_{(1)(2,3)}$ having analogous descriptions.

Points on $D_{(1,2)(3)}$ are stabilized by the $2$-element subgroup of
$\Ss_3$ generated by the transposition $\tau\,{=}\,(1,2)$: For a
generic point on $D_{(1,2)(3)}$, $\tau$ fixes the point and thus the
generating line. For the single point in
$D_{(1,2)(3)}\,{\cap}\,D_{\{0\}}$, $\tau$ fixes~$0$ and the line
$\langle (1,1,-2)\rangle$ pointwise. Analogously, we see that points
on $D_{(1,3)(2)}$ and on $D_{(1)(2,3)}$ are stabilized by the
transpositions $(1,3)$ and $(2,3)$, respectively.

On $D_{\{0\}}$, however, we find points whose stabilizers the nested
set stratification does not distinguish: Stabilizers for points on
$D_{\{0\}}$ are trivial except for those points on the intersections
with one of the other three divisors, {\em and\/} for $3$ additional
points
\[
   \psi_{12}=(0,\langle (1,-1,0)\rangle)\, \quad
   \psi_{13}=(0,\langle (1,0,-1)\rangle)\, \quad
   \psi_{23}=(0,\langle (0,1,-1)\rangle)\,
\]
The $\psi_{ij}$ are stabilized by transpositions $(i,j)$,
$1\,{\leq}\,i\,{<}\,j\,{\leq}\, 3$, respectively, since the
transpositions fix~$0$ and flip the lines in the second
coordinate. In fact, the transposition $(i,j)$,
$1\,{\leq}\,i\,{<}\,j\,{\leq}\,3$, acts on the open M\"obius band
$Y_{\Pi_3}$ like a ``central symmetry'' with fixed
point~$\psi_{ij}$.

\begin{center}
  \begin{picture}(0,0)%
    \includegraphics{cheburashkapoints.pstex}%
  \end{picture}%
  \input{cheburashkapoints.pstex_t}%
  \\
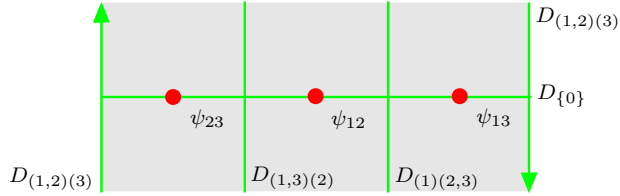
Figure 1. The nested set stratification $(Y_{\Pi_3}, \mfD)$.
\end{center}

\medskip

We provide here a glance on the already more complicated situation 
for $n\,{=}\,4$. Our picture below shows the stratification of the 
exceptional divisor~$D_{\{0\}}$, a real projective space of dimension~$2$,
as it emerges from the first blowup step in the De~Concini-Procesi 
construction, Bl$_{\{0\}}V$. 

\begin{center}
  \begin{picture}(0,0)%
    \includegraphics{cheburashkalines.pstex}%
  \end{picture}%
  \input{cheburashkalines.pstex_t}%
  \nopagebreak[4]\\
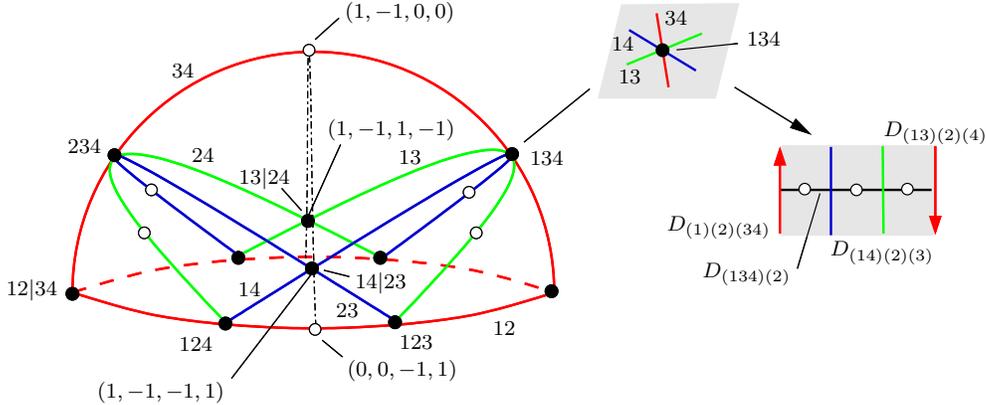
Figure 2. The stratification of $D_{\{0\}}$ after blowup of $\{0\}$ in~$V$.
\end{center}

We choose to place the intersection of $D_{\{0\}}$ with the hyperplane 
$H_{1,2}$ on the equator of the upper hemisphere model, and
thus obtain the stratification of $D_{\{0\}}$ by the braid arrangement
as depicted above. The double, respectively, triple intersections of 
hyperplanes in~$D_{\{0\}}$, e.g., $H_{1,2}\,{\cap}\,H_{3,4}$, respectively,
$H_{1,3}\,{\cap}\,H_{1,4}\,{\cap}\,H_{3,4}$, remain to be 
blown up in later steps, for triple intersections locally producing the 
situation that we studied above for $n\,{=}\,3$. 

We mark some points and lines on open strata that ought to be
distinguished by a stabilizer distinguishing stratification: For instance, 
the point on $D_{\{0\}}$ given by the line that is generated by the vector
$(0,0,-1,1)$ in $H_{1,2}$ should be distinguished from the open
stratum corresponding to $H_{1,2}$, since not only the transposition
$\tau=(1,2)$ but also $\sigma=(3,4)$ stabilizes this line. The same
goes for the (dashed) line obtained on $D_{\{0\}}$ as the intersection
with the plane spanned by the vectors $(1,-1,0,0)$ and $(0,0,-1,1)$.


\section{The nested set stratification of arrangement models}
\label{sect_nstrat}

\subsection{Points in $Y_{\GG}$} \label{ssect_pntYG}

Let $\AA$ be an arrangement of subspaces in a real vector
space~$V$, $\LL(\AA)$ its intersection lattice and
$\GG\,{=}\,\LL(\AA)\,{\setminus}\,\{\hat 0\}$ the maximal building 
set for~$\AA$.
We will encode points in the arrangement model~$Y_{\GG}$ into tuples
of points and lines in~$V$, a description that will prove to be
favorable for technical purposes.

A point $\omega$ in $Y_{\GG}$ will be written as
\begin{equation}\label{eqn_pntYG}
\omega \,\, = \,\, (x,H_1, l_1, H_2, l_2, \ldots, H_t,l_t)\,,
\end{equation}
where $x$ is a point in~$V$, the $H_i$ are elements in
$\GG\,{=}\,\LL\,{\setminus}\,\{\hat 0\}$, and the $l_i$ are lines in $V$.  
The point $x$ is
the first coordinate of $\omega$ when written as an element in the
product space on the right hand side of~(\ref{eqn_defmapYG}).
$H_1$ is the maximal lattice element that, as a subspace of~$V$,
contains~$x$. The line $l_1$ is orthogonal to $H_1$ and
corresponds to the coordinate entry of $\omega$ indexed by $H_1$
in $\P(V/H_1)$.  The lattice element $H_2$, in turn, is the
maximal lattice element that contains both $H_1$ and $l_1$. The
specification of lines $l_i$, i.e., lines that correspond to
coordinates of $\omega$ in $\P(V/H_i)$, and the construction of
lattice elements $H_{i+1}$, continues analogously for $i\geq 2$
until a last line $l_t$ is reached whose span with $H_t$ is not
contained in any lattice element other than the full ambient
space~$V$. Note, that if $H_t$ is a~hyperplane, then the line
$l_t$ is uniquely determined. The whole space $V$ can be
thought of as $H_{t+1}$.

Observe that the lattice elements $H_i$ are determined by the point
and the sequence of lines; we still choose to include the $H_i$ in
order to keep the notation more transparent.

To see that the description \eqref{eqn_pntYG} of a point $\omega$
in the arrangement model $Y_{\GG}$ is sufficient, we need to see
that the rest of the coordinates can be read off uniquely from the
coordinates $x,l_1,\dots,l_t$. The reconstruction can be
explicitly done as follows. Fixing $H_0:={0}$ and $l_0:=\langle
x\rangle$, the first coordinate of $\omega$ is $x$, and the
coordinate of $\omega$ indexed with $H\,{\in}\,\GG$, $\omega_H$, 
can be read from (\ref{eqn_pntYG}) as
\begin{equation} \label{eqn_expo}
  \omega_H\, \, = \, \, \langle l_j,H\rangle / H  \, \, \in \, \, \P(V/H)\, ,
\end{equation}
where $j$ is chosen from the index set $\{1,\ldots, t\}$ such that
$H\leq H_j$, but $H\not \leq H_{j+1}$.

To prove \eqref{eqn_expo} we need the following technical lemma.

\begin{lm}\label{lm4.1}
  Let $V$ be a vector space and $\wti H$, $H$ vector subspaces of $V$,
  such that $\wti H\subseteq H$. Let furthermore $(x_i)_{i=1}^\infty$
  be a~sequence of points in $V\setminus H$ such that the limit
  $\lim_{i\ra\infty}\langle x_i,\wti H\rangle=\Sigma$ exists in the
  corresponding Grassmannian.

  Assume that $\Sigma\not\subseteq H$, then $\lim_{i\ra\infty}\langle
  x_i,H\rangle=\langle\Sigma,H\rangle$; again the limit is understood
  with respect to the topology of the appropriate Grassmannian.
\end{lm}

\begin{pf}
 Let us split $V$ into the direct sum of linear subspaces:
\[V=\wti H\oplus(\wti H^\perp\cap H)\oplus H^\perp,\]
where $\wti H^\perp$, resp.\ $H^\perp$, denotes the orthogonal
complement of $\wti H$, resp.\ of~$H$.

Since $x_i\not\in\wti H$, we have $\dim\langle x_i,\wti H\rangle
=\dim\wti H+1$, hence $\dim\Sigma=\dim\wti H+1$, and therefore
there exists $v\in\wti H^\perp$, $v\neq 0$, such that
$\Sigma=\langle\wti H,v\rangle$.

Writing $x_i=a_i+b_i+c_i$, where $a_i\in\wti H$, $b_i\in\wti
H^\perp\cap H$, and $c_i\in H^\perp$, for all~$i$, we have
\begin{equation}\label{eqlm4.1}
\langle x_i,\wti H\rangle=\langle b_i+c_i,\wti H\rangle.
\end{equation}
Note that $b_i+c_i\in\wti H^\perp$, and $b_i+c_i\neq 0$. We can
scale $x_i$, such that $|b_i+c_i|=1$, and, after scaling $v$ and
changing $x_i$ to $-x_i$ for some appropriately chosen~$i$, we get
that $\lim_{i\ra\infty}(b_i+c_i)=v$. Denote
$\lim_{i\ra\infty}b_i=v_1$ and $\lim_{i\ra\infty}c_i=v_2$; these
limits exist since $b_i$ and $c_i$ are chosen in mutually
orthogonal linear subspaces. We certainly have
$\lim_{i\ra\infty}(b_i+c_i)=
\lim_{i\ra\infty}b_i+\lim_{i\ra\infty}c_i=v_1+v_2$, and
$v_1\in\wti H^\perp\cap H$, $v_2\in H^\perp$. Since $v\not\in H$,
we have $v_2\neq 0$, hence, for large~$i$, $|c_i|\geq|v_2|/2>0$.

We finish the proof by writing down two sequences of identities.
First,
\[\langle\Sigma,H\rangle=\langle\wti H,v,H\rangle=\langle v,H\rangle=
\langle v_1+v_2,H\rangle=\langle v_2,H\rangle,\]
where the second equality follows from $\wti H\subseteq H$, and the
fourth equality follows from $v_1\in H$. Second,
\[\lim_{i\ra\infty}\langle x_i,H\rangle=
\lim_{i\ra\infty}\langle c_i,H\rangle=
\langle\lim_{i\ra\infty}c_i,H\rangle=\langle v_2,H\rangle,\]
where the first equality follows from \eqref{eqlm4.1} and the fact
that $b_i\in H$. The second equality is the most interesting one, it
follows from the fact that the points $c_i$ lie in $H^\perp$,
and that the projectivization
map $\gamma:H^\perp\setminus\{0\}\ra\P(H^\perp)$, mapping a~point to
the line which it spans, is continuous.
\end{pf}

{\noindent {\bf Proof of \eqref{eqn_expo}.}} Choose a sequence
$(x_i)_{i=1}^\infty$, $x_i\in\MM(\AA)$, such that
$\lim_{i\ra\infty}\Phi(x_i)=w$ in $V\, \times \, \prod_{G\in\GG}\,
\P(V/G)$. This translates into
\[
\begin{cases}
  x=\lim_{i\ra\infty} x_i,\\
  \omega_G=\lim_{i\ra\infty}\Phi_G(x_i)=\lim_{i\ra\infty}\langle
  x,G\rangle/G.
\end{cases}
\]
Let us choose $H\,{\in}\,\GG$, and $j\in\{1,\ldots,t\}$, such that
$H\leq H_j$, but $H\not \leq H_{j+1}$. The identity \eqref{eqn_expo}
follows now from the following computation:
\[\lim_{i\ra\infty}\langle x_i,H\rangle=
\lim_{i\ra\infty}\langle \langle x_i,H_j\rangle ,H\rangle =
\lim_{i\ra\infty}\langle\langle l_j,H_j\rangle,H\rangle=
\lim_{i\ra\infty}\langle l_j,H\rangle,\] where the first and the
third equality are consequences of $H_j\subseteq H$, while the
second one follows from Lemma~\ref{lm4.1}. {\hfill
$\Box$\vspace{0.3cm}}

\subsection{Stabilizers of points in $Y_{\GG}$} \label{ssect_stabYG}

We now assume that our subspace arrangement carries the action of 
a~finite group $\Gamma$. As we discussed above, the
action extends to the arrangement model $Y_{\GG}$. When
considering stabilizers of the various actions we will include
indices into the notation that indicate the set on which the full
group is acting, e.g., we will write $\stab_V(y)$,
$\stab_{Y_{\GG}}(y)$ for the stabilizers of $y$ with respect to
the $\Gamma$-actions on~$V$ and on $Y_{\GG}$, respectively.

We take up the encoding of points in $Y_{\GG}$ from subsection
\ref{ssect_pntYG}, and derive a description for the stabilizer of
a point in $Y_{\GG}$:

\begin{prop} \label{prop_stab}
  Let an arrangement model $Y_{\GG}$ be equipped with a group
  action stemming from the action of a finite group $\Gamma$
  on the arrangement. Then for stabilizers of points
  $\omega\,{=}\,(x,H_1, l_1, H_2, l_2, \ldots, H_t,l_t)$ in
  $Y_{\GG}$ the following description holds:
\begin{equation}\label{eq_stab}
\stab_{Y_{\GG}} (\omega) \, \, = \, \, \stab_V(x)\, \cap \,
                               \stab_V(l_1) \, \cap \, \ldots
                               \, \cap \, \stab_V(l_t)\, ,
\end{equation}
where $\stab_V(l_i)$, $i\,{=}\,1,\dots, t$, denotes the subgroup
of elements $\gamma\,{\in}\,\Gamma$ with $\gamma(l_i)=l_i$, i.e., elements
preserving $l_i$ without necessarily fixing the line pointwise.
\end{prop}

\begin{pf}
  Using the description of points in $Y_{\GG}$ given in 
  subsection \ref{ssect_pntYG}, and the definition of the group
  action, we can describe the stabilizer of a~point
  $\omega\,{\in}\,Y_{\GG}$ as follows:
\begin{equation} \label{eqn_stab}
\stab_{Y_{\GG}} (\omega) \, \, = \, \, \stab_V(x)\, \cap \,
                               \stab_{\P(V/H_1)}(l_1) \, \cap \, \ldots
                               \, \cap \, \stab_{\P(V/H_t)}(l_t)\, ,
\end{equation}
where $\stab_{\P(V/H_i)}(l_i)$, $i\,{=}\,1,\dots, t$, translating
from the projective to the original lin\-ear setting, means
elements $\gamma\,{\in}\,\Gamma$ under which {\it both} $H_i$ and $l_i$ are
invariant:
\[
\stab_{\P(V/H_i)}(l_i)\,\, :=
\,\,\stab_V(H_i)\,\cap\,\stab_V(l_i)\, .
\]
Again, $\stab_V(H_i)$ denotes group elements that preserve $H_i$ but
do not necessarily fix~$H_i$ pointwise.

We show that
\begin{eqnarray*}
\stab_V(x) & \subseteq & \stab_V(H_1)\, , \quad \mbox{and } \\
\stab_V(H_i)\,\cap\, \stab_V(l_{i}) & \subseteq &
\stab_V(H_{i+1})\, , \quad \mbox{for }\, i=1,\ldots, t{-}1\, ,
\end{eqnarray*}
which, successively applied for $i\,{=}\,t{-}1$, $i\,{=}\,t{-}2$,
etc., reduces the right hand side of (\ref{eqn_stab}) to the right
hand side of \eqref{eq_stab}, since $A\cap B=A$, for any two sets $A$
and $B$, such that $A\subseteq B$.

For $\gamma\in \stab_V(x)$, $x$ in contained in $\gamma(H_1)\cap H_1$. But
$H_1 \supseteq \gamma(H_1)\cap H_1$ is assumed to be maximal in
$\GG\,{=}\,\LL\,{\setminus}\,\{\hat 0\}$ containing~$x$, thus, it 
follows from the fact that
$\GG$ is closed under taking intersections, that $\gamma(H_1)= H_1$.
Similarly for $\gamma\in \stab_V(H_i)\, \cap\,\stab_V(l_i)$:
$\gamma(H_{i+1})\cap H_{i+1}$ contains both $H_i$ and $l_i$, but
$H_{i+1}$ should be maximal in $\GG=\LL\,{\setminus}\,\{\hat 0\}$ 
with this property, hence $\gamma(H_{i+1})=H_{i+1}$.

Note additionally, that if $H_t$ is a~hyperplane, then
$\stab_V(H_t)=\stab_V(l_t)$, hence, in this case, $\stab_V(l_t)$ can be removed
from the right hand side of \eqref{eq_stab} without changing the
expression.
\end{pf}

\subsection{The divisors $D_G$, $G\,{\in}\,\GG$} \label{ssect_pntDG}

Recall from Section~\ref{sect_DPmodels} that the nested set stratification
$(Y_{\GG}, \mfD)$ on an arrangement model~$Y_{\GG}$ is given by
irreducible components of divisors and their intersections.
Our objective is to provide, in our special setting, a description of the
divisors $D_G$, $G\in \GG$,
that enables us to tell for a given point in the
arrangement model on which of these divisors it lies.

De~Concini~\& Procesi give a description of the divisors in terms
of affine and projective arrangement models for ``smaller''
arrangements. To keep track of the respective settings, we provide
arrangement models with an additional index that specifies the
ambient space of the original arrangement, and we indicate
projective models by a bar, e.g., in presence of other arrangement
models we will now write $Y_{V,\GG}$ for the affine and
$\ol{Y}_{V,\GG}$ for the projective model of the previously
considered arrangement.

In our special setting the description of divisors by De~Concini
\& Procesi reads as follows:

\begin{prop}
  \cite[Thm.\,4.3, Rem.\,4.3.(1)]{DP} Let $\AA$ be an essential
  arrangement of subspaces, $\GG$ the maximal building set,
  $\GG\,{=}\,\LL(\AA)\,{\setminus}\,\{\hat 0\}$, and $Y_{V,\GG}$ 
  the corresponding arrangement
  model. For the irreducible divisors $D_G$, $G\in \GG$, there are
  natural isomorphisms:
\begin{eqnarray}
D_{\{0\}} & \cong  & \ol{Y}_{V,\GG} \,, \label{eqn_D0}\\
D_{\{G\}} & \cong  & \ol{Y}_{V/G,\GG_{\leq G}}
                                     \, \, \times \, \,
Y_{G,\GG_{>G}}, \qquad \mbox{for }\,G\not = \{0\}\, .
\label{eqn_DG}
\end{eqnarray}

Here, $\ol{Y}_{V/G,\GG_{\leq G}}$ is the projective model for the
quotient arrangement $\AA/G:=\{H/G\, | \, H\,{\in}\,\AA,
H\,{\supseteq}\,G\}$ with (maximal) building set $\GG_{\leq
  G}=\{H\,{\in}\,\GG\,|\, H \leq G\}$, and $Y_{G,\GG_{>G}}$ is the
affine model for the restricted arrangement $\AA\cap
G:=\{H\,{\cap}\,G\, | \, H\,{\in}\,\AA\}$ with (maximal) building set
$\GG_{>G}=\{H\,{\in}\,\GG\,|\, H\,{>}\,G\}$.
\end{prop}

The projective model $\ol{Y}_{V,\GG}$, in fact, is isomorphic to the
inverse image of $\{0\}$ when projecting $Y_{V,\GG}$ to $V$, the first
coordinate of its ambient space~\cite[Thm.4.1]{DP}. Hence,
$\omega\,{\in}\, D_{\{0\}}$ if and only if $\omega_{\{0\}}\,{=}\,0$,
in other words
\begin{equation} \label{eqn_charD0}
       \omega \in  D_{\{0\}} \, \, \Leftrightarrow \, \,
       \omega \in  Y_{V,\GG} \,\cap \,\left(
               \{0\} \,\times\,\prod_{G\in\GG}\, \P(V/G)\,\right).
\end{equation}
It is a description of this type that we want to achieve for the other
divisors, $D_{G}$, $G\,{\not=}\,\{0\}$, as well.

To this end, note that the right hand side of (\ref{eqn_DG}) can
be considered as a subspace of
\[
       \{0\}\, \times \, \prod_{H\in \GG_{\leq G}} \, \P(V/G \big/ H/G)
           \,\, \times \,\, G \, \, \times \,
                        \prod_{H\in \GG_{>G}} \, \P(G/H)\,.
\]
For $K\,{\in}\,\GG_{>G}$, we can ``expand'' the factor $\P(G/K)$
by a diagonal map
\[
\P(G/K)\,\,\lra  \prod_{H\in \GG \atop H \vee G=K } \,
\P(G/(H\,{\vee}\,G)) \, ,
\]
and thus interpret  $D_G$ as a subset of
\[
   U_G\, \, := \, \, G\, \, \times
               \, \prod_{H\in \GG_{\not\leq G}} \, \P(G/ (H\,{\vee}\,G))
           \,\, \times \,\,
                  \prod_{H\in \GG_{\leq G}} \, \P(V/H) \, .
\]
With $G/(H\,{\vee}\,G)\,{\cong}\, \langle G,H\rangle/H$, $U_G$
can be considered a subspace of the ambient space
$V\, \times\, \prod_{H\in \GG} \P(V/G)$ of the arrangement model.

We thus can state our description of divisors $D_G$:

\begin{prop} \label{prop_charDG}
  Let $\AA$ be an essential arrangement of subspaces, $\GG$
  the maximal building set, $\GG\,{=}\,\LL(\AA)\,{\setminus}\,\{\hat 0\}$, 
  and $Y_{\GG}$ the
  corresponding arrangement model. The irreducible divisors $D_G$, $G
  \in \GG$, are intersections of $Y_{\GG}$ with the product spaces
  $U_G$, where the $U_G$ are obtained by restricting those factors of
  the original ambient space of $Y_G$ which are indexed with
  $H\,{\in}\,\GG_{\not\leq G}$:
\begin{eqnarray*}
D_{G} & = & Y_{\GG} \, \, \cap \, \, U_G \\
      & = & Y_{\GG} \, \, \cap \, \, \left(
               G\, \, \times
      \,\prod_{H\in \GG_{\not\leq G}} \, \P(\langle G,H\rangle/H)
           \,\, \times \,\,
                  \prod_{H\in \GG_{\leq G}} \, \P(V/H) \right)  \, .
\end{eqnarray*}
\end{prop}

\begin{pf}
Observe first that the description for $D_{\{0\}}$ given in
(\ref{eqn_charD0}) coincides with the one stated in the
Proposition: intersecting $Y_{\GG}$ with $U_{\{0\}}$ restricts the
first coordinate to~$0$.

For $G\,{\not=}\,\{0\}$, we start with the description of $D_G$ in
(\ref{eqn_DG}) and see from the reasoning above that any element
in $D_G$ is contained in $U_G$. For the converse, let
$\omega\,{=}\,(x,H_1, l_1, H_2, l_2, \ldots, H_t,l_t)$ be
contained in $Y_{\GG}\,{\cap}\,U_G$. From $\omega\,{\in}\, U_G$ we
conclude that $x\,{\in}\,G$, hence $H_1\,{\geq}\,G$. Assuming for
the moment that $H_1\,{\gneq}\,G$, we look at the component of
$\omega$ indexed by $H_1$. Using the expansion of $\omega$ from
(\ref{eqn_expo}) and the fact that $\omega\,{\in}\, U_G$, we see
that
\[
   \omega_{H_1}\, \, =\, \, \langle l_1, H_1 \rangle / H_1
               \, \, \in \, \, \P(G / H_1)\, ,
\]
hence $l_1\,{\subseteq}\,G$. This implies that $H_2$ is larger or
equal $G$, for, if it were not, $H_2\,{\vee}\, G\,{\gneq}\, H_2$ would
contain both $H_1$ and $l_1$ in contradiction to $H_2$ being maximal
with this property.

We conclude that there is an index $k\,{\in}\,\{1,\ldots,t\}$ with
$H_k\,{=}\,G$, and can thus split the point/lines description of
$\omega$ into
\[
\omega\,\, = \, \,\Big(\, (x,H_1, l_1, H_2, l_2, \ldots,l_{k-1},G),
                   (l_k, H_{k+1}, \ldots, H_t,l_t)\, \Big)\, .
\]
The first tuple clearly describes an element in $Y_{G,\GG_{>G}}$.
We rewrite the second tuple as follows:
\[
   (0_{V/G}, l_k, H_{k+1}/G, \ldots, H_t/G,l_t)       \, .
\]
With $l_j$ being orthogonal to $G$, hence $l_j\in \P(V/G)$, we can
then interpret it as an element of $\ol{Y}_{V/G, \GG_{\leq G}}$.
With (\ref{eqn_DG}) we thus conclude that $\omega\in D_{G}$.
\end{pf}

\subsection{Open strata of the nested set stratification}
\label{ssect_nsstr}

We will provide a characterization of points on open strata of the
nested set stratification of~$Y_G$ in terms of their point/line
encoding described in subsection~\ref{ssect_pntYG}.

To fix some notation, let us denote by
$D_{G_1,\ldots,G_m}^{\circ}$ the open stratum in $(Y_{\GG}, \mfD)$
that lies in the intersection of divisors $D_{G_1}, \ldots,
D_{G_m}$, but on no other divisors indexed with building set
elements. Recall that the index set $\{G_1,\ldots,G_m\}$ is
$\GG$-nested, which in our context, i.e., for the maximal building
set, means that it is a chain in~$\LL(\AA)$. We tacitly assume
that the $G_i$ are listed in a descending order:
$G_1\,{>}\,\ldots\,{>}\,G_m$.

\begin{prop} \label{prop_nsstr}
Let $Y_{\GG}$ be an arrangement model with nested set
stratification $\mfD$. A point $\omega\,{\in}\,Y_{\GG}$ is
contained in the open stratum of $\mfD$ indexed with the nested
set $\TT\,{=}\,\{G_1,\ldots, G_m\}$ if and only if the spaces in
$\TT$ coincide with the spaces occurring in the point/line
description of~$\omega$:
\begin{eqnarray*}
       \omega \in D_{G_1,\ldots,G_m}^{\circ} &
                               \Longleftrightarrow &
       \omega=(x,G_1,l_1,\ldots, G_m,l_m) \, ,
\end{eqnarray*}
where on the right hand side the usual restrictions for coordinates of
a point/line tuple as in\/~{\rm (\ref{eqn_pntYG})} apply.
\end{prop}

\begin{pf}
First observe that the claim holds for points~$\omega$ in the big
open stratum $Y_{\GG}\,{\setminus}\,D\,{=}\,\MM(\AA)$, that is for
$m=0$ : The indexing nested set is empty, and the point/line
description for~$\omega$ reduces to the point entry
$x\,{\in}\,\MM(\AA)$.

We can thus assume that $\omega\,{\in}\,D$, in particular,
$\omega$ is contained in some open stratum in~$D$, say
\[
    \omega\, \, \in \, \, D_{G_1,\ldots,G_m}^{\circ},
\]
where we remind that the $G_i$ are indexed in descending order,
and $m\geq 1$.

At the same time, $\omega$ has a point/line description, say
\[
     \omega\, \, = \, \,(x,H_1,l_1,\ldots,H_t,l_t)\, ,
\]
where $H_1,\ldots,H_t\,{\in}\,\GG$, $x\,{\in}\,H_1$, and
$l_i\,{\in}\,\P(V/H_i)$, for $i\,{=}\,1,\ldots,t$. We show in the
following that the descending chains $G_1\,{>}\,\ldots\,{>}\,G_m$
and $H_1\,{>}\,\ldots\,{>}\,H_t$ coincide, in particular implying
$m=t$.

\noindent
{\bf Step 1:} {\em The maximal elements of the chains coincide:
$H_1\,{=}\, G_1$.}
\\
With $\omega\,{\in}\,D_{G_1}$, we know by
Proposition~\ref{prop_charDG} that $x\,{\in}\, G_1$; but $H_1$ is
maximal with this property, hence, $H_1\,{\geq}\,G_1$.

We want to see, that $\omega\,{\in}\,D_{H_1}$. Using again
Proposition~\ref{prop_charDG} and the expansion of~$\omega$
in~(\ref{eqn_expo}), we have to check that $x\,{\in}\,H_1$, and
that for any $H\,{\not \leq}\,H_1$ the coordinate
$\omega_H\,{=}\,\langle \langle x\rangle,H\rangle/H$ is a point in
$\P(\langle H_1,H\rangle/H)$. With $\langle
x\rangle\,{\subseteq}\,H_1$ this is obviously the case.

We conclude that $H_1\,{\in}\,\TT$, hence, $H_1\,{\leq}\,G_1$ by
maximality of $G_1$ in~$\TT$. This yields our claim. In
particular, we see that $t\geq 1$.

\noindent
{\bf Step 2:} {\em Assume $H_j\,{=}\,G_j$ for $j\,{=}\,1,\ldots, i$,
  and $i\,{\lneq}\,t$. Then $m\,{\geq}\,i{+}1$ and
  $H_{i+1}\,{=}\,G_{i+1}$.}
\\
Here, we first want to see, that $\omega\,{\in}\,D_{H_{i+1}}$. For
this we need to check that $x\,{\in}\, H_{i+1}$, and that for any
$H\,{\not \leq}\,H_{i+1}$ the coordinate $\omega_H\,{=}\,\langle
l_j,H\rangle/H$ is a point in $\P(\langle H_{i+1},H\rangle/H)$.
The line $l_j$ depends on $H$ (compare (\ref{eqn_expo})), but for
any $H$ in question its index $j$ is strictly less
than~$i\,{+}\,1$. From the point/line description for $\omega$ we
see that $x\,{\in}\,H_1\,{\subseteq}\, H_{i+1}$. With
$l_j\,{\subset}\, H_{j+1}\,{\subseteq}\,H_{i+1}$ we conclude that
$\langle l_j,H\rangle/H\,{\in}\,\P(\langle H_{i+1},H\rangle/H)$,
hence $\omega\,{\in}\,D_{H_{i+1}}$.

Since $H_{i+1}$ belongs to the nested set~$\TT$,
$H_{i+1}\,{<}\,H_i\,{=}\,G_i$, implies that, in fact,
$m\,{\geq}\,i{+}1$ and $H_{i+1}\,{\leq}\,G_{i+1}$.

To obtain equality we write out the condition on the coordinate of
$\omega$ indexed with~$H_i$ that results from
$\omega\,{\in}\,D_{G_{i+1}}$: $\omega_{H_i}\,{=}\, \langle
l_i,H_i\rangle/H_i\,{\in}\, \P(\langle
G_{i+1},H_i\rangle/H_i)=\P(G_{i+1}/H_i)$.

We conclude that $l_i\,{\subseteq}\, G_{i+1}$. Moreover,
$G_i\,{\subseteq}\,G_{i+1}$ by descending order on $\TT$. But
$H_{i+1}$ is maximal in $\GG$ containing both $H_i\,{=}\,G_i$ and
$l_i$, hence $H_{i+1}\,{\geq}\,G_{i+1}$, from which our claim
follows.

\noindent
{\bf Step 3:} {\em $m\,{=}\,t$, and hence the chains coincide.}
\\
From Steps (1) and (2) we conclude that $m\,{\geq}\,t$. Let us
assume that $m\,{>}\,t$, in particular,
$\omega\,{\in}\,D_{G_{t+1}}$. We conclude from the resulting
condition on the coordinate indexed by $H_t$,
$\omega_{H_t}\,{=}\,\langle l_t,H_t\rangle/H_t\,{\in}\,\P(\langle
G_{t+1},H_t\rangle/H_t)=\P(G_{i+1}/H_i)$, that both $l_t$ and
$H_t\,{=}\,G_t$ are contained in $G_{t+1}$ which contradicts the
fact that the point/line description of $\omega$ was terminated
after the $t$-th step. Hence $m\,{=}\,t$, and the chains
$G_1\,{>}\,\ldots\,{>}\,G_t$ and $H_1\,{>}\,\ldots\,{>}\,H_t$
coincide.
\end{pf}


\section{A stabilizer distinguishing stratification of $Y_{\Pi_n}$}
\label{sect_stdist}

\subsection{Adding strata} \label{ssect_addstr}
On our way to construct a stabilizer distinguishing stratification for
$\YP$ we first analyze the locus of {\em lines\/} in $\R^n$ that are
stabilized by a given element in~$\Ss_n$.  Let $\pi\,{\in}\,\Ss_n$,
and, restricting the permutation action, consider $\R^n$ as
a~representation space of the cyclic group $\langle \pi \rangle$.
In  $\R^n$ we have, on one hand, the linear subspace
$T_1(\pi)\,{=}\,{\rm Fix}(\pi)$, the locus of lines that are pointwise
fixed by~$\pi$, on the other hand, we have the subspace $T_{-1}(\pi)$,
the locus of lines that are flipped by~$\pi$.
We can characterize lines in $\R^n$ that
are invariant under $\pi\,{\in}\, \Ss_n$ as follows:

\begin{prop} \label{propo_stabl}
Let $\pi\,{\in}\,\Ss_n$ and $S(\pi):=T_1(\pi)\,{\cup}\,T_{-1}(\pi)$.
For a given line $l$ in $\R^n$,
\[
\pi \in \stab\,(l) \quad \Longleftrightarrow \quad
l\,{\subseteq}\, S(\pi)\, .
\]
\end{prop}

We would like to emphasize that $S(\pi)$ is defined as a union of
$T_1(\pi)$ and $T_{-1}(\pi)$, {\it not} as their span.

Let us now describe stratifications of the orthogonal complements
$G^{\perp}$ of subspaces~$G$ in~$\Pi_n$. For such $G$, and for any
$\pi\,{\in}\, \Ss_n$, define $S(\pi,G)\,{:=}\,S(\pi)\,{\cap}\,
G^{\perp}$. Then,
\[
        \mfS_G\, \, := \, \, \big\{ S(\pi,G) \big\}_{\pi\in \Ss_n}
\]
is a stratification of $G^{\perp}$. Unlike the restriction of the
braid arrangement stratification to $G^{\perp}$, it distinguishes
stabilizers of points as well as stabilizers of lines.

We propose a construction for subsets in  real
arrangement models $Y_{\GG}$ that takes unions of
linear subspaces in~$\R^n$ as input data. It is inspired
by the description of divisors $D_G$, $G\,{\in}\,\GG$, that we
presented in Proposition~\ref{prop_charDG}. Taking spaces
$S(\pi,G)\,{\times}\,G$, $G\,{\in}\,\GG$, $\pi\,{\in}\,S_n$, with
$S(\pi,G)$ as defined above, our construction will
provide us with the additional maximal strata in $\YP$ for obtaining a
stabilizer distinguishing stratification.

\begin{df} \label{def_addstr}
Let $Y_{V,\GG}$ be an arrangement model, and $W\,{=}\,\{W_1,\ldots,W_m\}$
a family of real linear subspaces in $V$. Define a subset $B(W)$ in $Y_{\GG}$
by
\begin{eqnarray*}
B(W) & := & Y_{\GG} \, \, \cap \, \, \left( \,
               \bigcup W\, \, \,\,\times \,
               \prod_{{H\in \GG, H\not \supseteq W_i}\atop
                      {\mbox{{\tiny for any} }\, W_i \in W}}
                 \, \P(\langle W,H\rangle/H)
           \,\, \times \,\,
                  \prod_{{H\in \GG, H \supseteq W_i}\atop
                      {\mbox{{\tiny for some }}\, W_i \in W}}
                   \, \P(V/H) \, \right)  \, ,
\end{eqnarray*}
where $\P(\langle W,H\rangle/H)$ stands for the projectivization of
$\bigcup_{i=1}^m \langle W_i,H\rangle/H$.
\end{df}

We now can refine the nested set stratification $\mfD$ of $\YP$ so as to
obtain a stabilizer distinguishing stratification. As before, we describe
the stratification by listing its maximal strata:
\begin{equation} \label{eqn_strat}
\mfB\, \, :=\, \,
        \left\{ \, \, \big( \, D_G \, \big)_{G\in \Pi_n},\,
   \big( \, B(\,S(\pi,G)\times G\,)\,  \big)_{{G\in \Pi_n}, \pi\in \Ss_n}
        \right\}\, ,
\end{equation}
where in the second family of strata we only consider those
with $\{0\}\,{\subsetneq}\,S(\pi,G)\subseteq G^{\perp}$.


\subsection{${\mathbf (Y_{\Pi_n}, \mfB)}$ is stabilizer distinguishing}
\label{ssect_Bstbd}

We can now state one of the main results of this article:

\begin{thm} \label{thm_Bstbd}
  The stratification $\mfB$ for the arrangement model $\YP$ defined
  in~(\ref{eqn_strat}) is stabilizer distinguishing, i.e., the
  stabilizer of a point $\omega\in \YP$ is completely determined by
  the open stratum of $\mfB$ that contains~$\omega$.
\end{thm}

\begin{pf}
We pick a point $\omega\,{=}\,(x,G_1,l_1,\ldots, G_t,l_t)$ in
$\YP$, and assume that we have the complete list of maximal strata
in $\mfB$ which contain~$\omega$. We want to show that the
stabilizer of $\omega$ is fully determined by this list.

Note first that by Proposition~\ref{prop_nsstr} our list of strata
contains the divisors $D_{G_1},\ldots, D_{G_t}$, and no other
divisors of this type. This means that we can read of from the
list the elements $G_1,\dots,G_t$ for the point/line description
of~$\omega$.

Assume $\omega\,{\in}\,B(S(\pi,G_i)\,{\times}\, G_i)$, for some
$G_i$, $i\,{\in}\,\{1,\ldots,t\}$. With
Definition~\ref{def_addstr}, and
$S(\pi,G_i)\,{\times}\,G_i\,{\supseteq}\,\,G_i$,  this puts the
following restriction on the coordinate of $\omega$ that is
indexed by~$G_i$:
\[
    \omega_{G_i} =  \langle l_i,G_i\rangle /G_i
                 \,\, \in \, \,
    \P(\langle S(\pi,G_i)\,{\times}\,G_i , G_i \rangle / G_i)\, .
\]
We conclude that $l_i\,{\subseteq}\,S(\pi,G_i)$, in particular, $\pi$
stabilizes $l_i$.

From the strata $B(S(\pi,G_i)\,{\times}\, G_i)$, that occur on our
list for a fixed space $G_i$, $i\,{\in}\,\{1,\ldots,t\}$, we can
read off a subset $\Gamma_i$ of $\stab(l_i)$. Namely, for each 
$i\,{\in}\,\{1,\ldots,t\}$, $\Gamma_i$
consists of all $\pi$ such that
$\omega\,{\in}\,B(S(\pi,G_i)\,{\times}\, G_i)$. 

Let us assume that, when constructing $\Gamma_i$ from our list of strata
for $\omega$, we actually missed some elements of~$\stab(l_i)$:
let $\sigma\,{\in}\, \stab(l_i){\setminus}\Gamma_i$. Then
$l_i\,{\subseteq}\,S(\sigma,G_i)$, but
$\omega\,{\not \in}\, B(S(\sigma,G_i)\,{\times}\, G_i)$. By definition of the
additional maximal strata we conclude that there exists a subspace
$H\,{\in}\,\Pi_n$, which does not contain any of the spaces in
$S(\sigma,G_i)\,\times\, G_i$, such that
\begin{equation} \label{eq_5.2.1}
   \omega_{H} = \langle l_j,H\rangle /H
                 \,\, \not \in  \, \,
    \P(\langle S(\sigma,G_i)\,{\times}\,G_i , H \rangle / H)\, .
\end{equation}
The line index $j$ depends on $H$, but in any case, $j>i$: for
$j<i$, $l_j\,{\subseteq}\, G_i$, and for $j=i$,
$l_i\,{\subseteq}\, S(\sigma,G_i)$, and the condition on
$\omega_H$ for $\omega$ being contained in
$B(\,S(\sigma,G_i)\,{\times}\, G_i\,)$ would be fulfilled.

It follows from \eqref{eq_5.2.1} that
$l_j\,{\not\subseteq}\,S(\sigma,G_i)$. Since $l_j$ is orthogonal
to~$G_i$, it implies $\sigma\,{\not\in}\, \stab(l_j)$, and, in
particular, $\sigma\,{\not\in}\, \bigcap_{i=1}^t\, \stab (l_i)$.
Hence, even if for some~$i$, $\Gamma_i\subsetneq\stab(l_i)$, once the full
intersection is taken, this is rectified:
\[
 \bigcap_{i=1}^t\, \Gamma_i \, \, = \, \,  \bigcap_{i=1}^t\, \stab (l_i)\, .
\]
With the description of $\stab(\omega)$ from
Proposition~\ref{prop_stab}, and $\stab(x)$ being determined by the
partition pattern of $x$, hence by $G_1$, we can conclude that the
list of strata in $\mfB$ containing $\omega$ actually determines the
stabilizers of~$\omega$.
\end{pf}


\subsection{$Y_{\Pi_3}$ revisited}

Let us have a look at the stratification $\mfB$ on $Y_{\Pi_3}$ and see how it
resolves the problem raised in~\ref{ssect_chebpnts}, namely to distinguish
stabilizers of points by means of a stratification.

To start with, we have to identify those spaces $S(\pi,G)\,{\times}\,G$
for $G\,{\in}\,\Pi_3$, $\pi\,{\in}\,\Ss_3$, that give raise to new strata
$B(S(\pi,G)\,{\times}\,G)$. We claim that the only interesting case occurs
for $\pi$ a transposition, $\pi=(i,j)$, $1{\leq}i{<}j{\leq}3$, and
$G\,{=}\,\{0\}$.

We have $S(\pi)\,{=}\,H_{i,j}\,{\cup}\,H_{i,j}^-$, where we denote
hyperplanes of $\AA_{n{-}1}$ in $V$ by $H_{i,j}$, just as for the original
(non-essential) arrangement in $\R^3$, and their orthogonal
complements by~$H_{i,j}^-$. With $S(\pi,\{0\})\,{=}\,S(\pi)$, we obtain
new strata
\[
    B_{(i,j)}\, \, = \, \, B(S((i,j),\{0\})\,{\times}\,\{0\}) \, \,= \, \,
    Y_{\Pi_3}\, \cap \, \big(   (H_{i,j} \cup H_{i,j}^-) \, \times \,
                               (\P(H_{i,j}) \cup \P(H_{i,j}-))
                        \big).
\]
In terms of the pointwise description for $Y_{\Pi_3}$ that we gave
in~\ref{ssect_chebpnts} this reads
\begin{eqnarray*}
    B_{(1,2)} &  =  & \{\,(x,\langle x\rangle)\,|\, x_1=x_2\neq 0
                            \mbox{ or }\,x_1=-x_2\neq 0 \, \}  \\
              &      & \, \qquad \cup \,\,
    \{\,(0,\langle(1,1,-2)\rangle),(0,\langle(1,-1,0)\rangle)\,\}\, ,
\end{eqnarray*}
analogously for $B_{(1,3)}$, $B_{(2,3)}$.
Hence, as opposed to the nested set stratification~$\mfD$, the stratification
$\mfB\,{=}\,\{(D_G)_{G\in \Pi_3}, B_{(1,2)},  B_{(1,3)},  B_{(2,3)}\}$
distinguishes the points $\psi_{i,j}$, $1\,{\leq}\,i\,{<}\,j\,\leq3$ from
the rest of the divisor~$D_{\{0\}}$.

\vspace{0.3cm}
\begin{center}
  \begin{picture}(0,0)%
    \includegraphics{mobius.pstex}%
  \end{picture}%
  \input{mobius.pstex_t}%
  \\
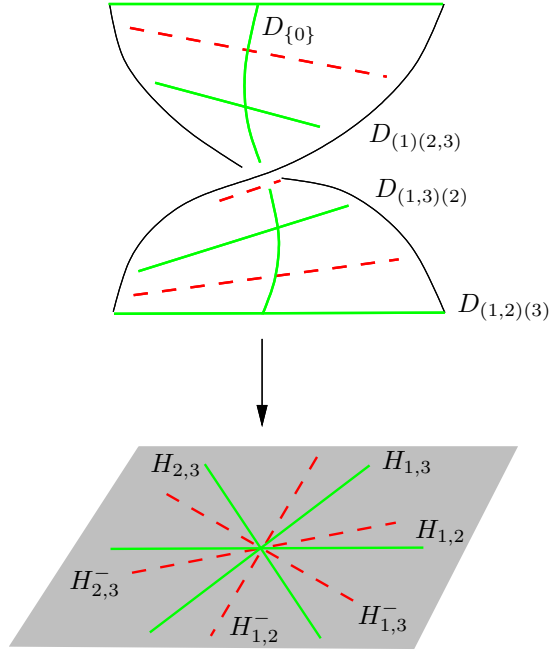
Figure 3. The stratification $(Y_{\Pi_3},\mfB)$.
\end{center}


\section{A combinatorial framework for describing stabilizers}
\label{sect_combfr}

In this section we develop a combinatorial framework for
describing stabilizers of points on the De~Concini-Procesi
arrangement model~$Y_{\Pi_n}$ with respect to the $\Ss_n$-action.
In Section~\ref{sect_stabYP} we will use this description to prove
that the stabilizers of points of $Y_{\Pi_n}$ are isomorphic to
direct products of ${\Z_2}$.

\subsection{Diagrams over families of cubes}

\begin{df} $\,$
\begin{enumerate}
\item Let $I$ be a finite, possibly empty set of positive integers. We
  call the collection of all subsets of $I$ (including the empty
  subset) an~{\bf $I$-cube}. Reversely, given an~$I$-cube $K$, we call
  $I$ the {\bf index set} of~$K$.
\item Let $t$ be a~positive integer. A~{\bf $t$-family of cubes} is a
  collection $\CC=\{K_1,\dots,K_p\}$, where, for each $j=1,\dots,p$,
  $K_j$ is an $I(j)$-cube, for some $I(j)\subseteq\{1,\dots,t\}$.
\end{enumerate}
\end{df}

One can make use of geometric intuition by thinking of an $I$-cube
as a coordinate $0/1$-cube with $I$ indexing the set of
``directions'' of the cube. The $\emptyset$-cube is simply the
point at the origin. For every $n\geq\max(I)$, the $I$-cube can be
imbedded as a coordinate $0/1$-cube in $\R^n$, and our object
is the equivalence class of all these imbeddings.

Let $K$ be an $I$-cube, to discriminate from other $I$-cubes, we
write elements of $K$ as pairs $(K,S)$, for $S\subseteq I$.  We
denote $\vrt(K)=\{(K,S)\,|\,S\subseteq I\}$, and refer to its
elements as vertices of $K$. When it is clear which cube we are
in, we may choose to skip $K$, and call $S$ itself a vertex
of~$K$.

Note also that a $t$-family of cubes is simply specified by a
function $I:[p]\rightarrow 2^{[t]}$, and that if $\tilde t>t$,
then every $t$-family of cubes is also a~$\tilde t$-family. For
$\CC=\{K_1,\dots,K_p\}$ we denote $\vrt(\CC)=\bigcup_{i=1}^p
\vrt(K_i)$, and refer to its elements as vertices of $\CC$.

\begin{df} \label{df6.2} $\,$
\begin{enumerate}
\item Let $\CC$ be a~$t$-family of cubes, $\CC=\{K_1,\dots,K_p\}$, and
  let $n$ be a positive integer. An~{\bf $n$-diagram} $\DD$ over $\CC$
  is a~partition of the set $[n]$ into $|\vrt(\CC)|$ blocks, some blocks 
  may be  empty, and an~assignment of the blocks of
  this partition to vertices of $\CC$; in other words, it is a~function
\begin{equation}\label{eq6.1}
\begin{array}{ccccc}
\DD&:&[n]&\longrightarrow& \vrt(\CC),\\
&&k&\mapsto& (K_{\alpha(k)},v_k),
\end{array}
\end{equation}
where $\alpha(k)\,{\in}\,[p]$ specifies the index of the cube
and $v_k\,{\subseteq}\,I(\alpha(k))$ the vertex of $K_{\alpha(k)}$
assigned to $k$.
\item For a vertex $(K,v)$ of $\CC$, we call the set $\DD^{-1}(K,v)$
  the {\bf fiber} of $\DD$ over $(K,v)$. For an $I$-cube $K$ in $\CC$,
  the fiber of $\DD$ over $K$ is defined as the union of the fibres of
  the vertices of $K$:
  $$\DD^{-1}(K):=\bigcup_{v\subseteq I}\DD^{-1}(K,v).$$
\end{enumerate}
\end{df}

\begin{center}
  \begin{picture}(0,0)%
    \includegraphics{cubes.pstex}%
  \end{picture}%
  \input{cubes.pstex_t}%
   \\
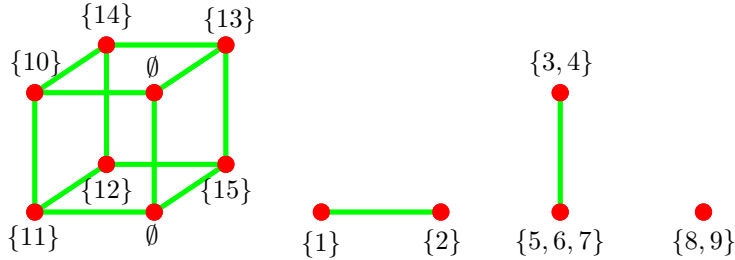
Figure 4. An example of a 15-diagram over a 3-family of cubes.
\end{center}

As yet another piece of notation, let $\rho(\DD)\vdash n$ be the set
partition with blocks being the fibers of $\DD$ over the vertices of
$\CC$, i.e., $\rho(\DD)=\{\DD^{-1}(K,v)\}_{(K,v)\in\vrt(\CC)}$, where
we disregard all the empty blocks in the set on the right hand side.

\subsection{Automorphism groups}

There is a standard $\Z_2^n$-action on an $[n]$-cube: it is generated
by reflections with respect to $n$ hyperplanes, which are parallel to
the facets of the cube, and which go through the center of the cube.
A~technically convenient way to describe this action is to think of
the vertices of an $[n]$-cube as vectors in an $n$-dimensional vector
space over the field $\F_2$, again denoted $\Z_2^n$, and the action as
parallel translations by vectors in $\Z_2^n$ (i.e., generated by
parallel translations with respect to the coordinate vectors).

For a subset $I\subseteq [n]$, let $\Z_2^{I}$ denote the
corresponding coordinate subspace of $\Z_2^n$, and let
$\proj_{I}:\Z_2^n\rightarrow \Z_2^{I}$ denote the projection onto
$\Z_2^{I}$ which simply "forgets" the coordinates with indices
outside of~$I$.

The following definition generalizes these actions to the case of
diagrams over families of cubes.

\begin{df} \label{df6.3}
  Let $\DD$ be an $n$-diagram over a~$t$-family of cubes
  $\CC=\{K_1,\dots,K_p\}$. We define the {\bf group of automorphisms}
  of $\DD$, which we denote $\aut(\DD)$, as follows:  $\aut(\DD)$
  consists of all permutations $\pi\in\Ss_n$, such that
\begin{enumerate}
\item[i)] $\pi_{\DD^{-1}(K_j)}\in\Ss_{\DD^{-1}(K_j)}$, for all
$j=1,\dots,p$, i.e., $\pi$ preserves the fibers over cubes;
\item[ii)] there exists (not necessarily unique) $\sigma\in\Z_2^t$,
  such that
\begin{equation}\label{eq6.2}
  v_{\pi(k)}=\sigma_{\alpha(k)}(v_k),\quad \text{ for all }\,\,
  k\in\{1,\dots,n\},
\end{equation}
where $\sigma_j=\proj_{I(j)}(\sigma)$, for all $j\in\{1,\dots,p\}$,
and where $v_k$ and $\alpha(k)$ are as in \eqref{eq6.1}. In other
words, $\pi$ maps fibers to fibers according to a~uniform scheme
obtained by restricting $\sigma$ to the cubes in the family $\CC$.
\end{enumerate}

\end{df}

\begin{rem} \label{rem6.4}
  Maps between fibers of an $n$-diagram $\DD$ over a $t$-family of
  cubes $\CC$, which are induced by an element $\pi\in\aut(\DD)$,
  must be bijections.
\end{rem}

Indeed, let $K$ be an $I$-cube in $\CC$, let $v\subseteq I$, and let
$\sigma\in\Z_2^t$ be associated to $\pi$ by Definition \ref{df6.3}
ii), then, by \eqref{eq6.2}, we have

$$\pi(\DD^{-1}(K,v))\subseteq\DD^{-1}(K,\proj_I(\sigma)(v)),$$
while
$$\pi(\DD^{-1}(K,\proj_I(\sigma)(v)))\subseteq
\DD^{-1}(K,\proj_I(\sigma)^2(v))=\DD^{-1}(K,v).$$

Since $\pi$ is injective, its restrictions are injective as well,
hence we can conclude that $\pi$ restricts to a bijection between
$\DD^{-1}(K,v)$ and $\DD^{-1}(K,\proj_I(\sigma)(v))$.

\begin{lm} \label{lm_stab} $\,$
\begin{enumerate}
\item For $x\in\R^n$, the stabilizer of $x$ under the $\Ss_n$-action is
  the Young subgroup of $\Ss_n$ indexed by the set partition of $[n]$,
  which is induced by the coordinates of $x$. One can represent this
  Young subgroup as an~automorphism group of an~$n$-diagram over
  a~0-family of cubes.
\item For a line $l\subseteq\R^n$, the stabilizer of $l$ under the
  $\Ss_n$-action can be represented as an~automorphism group of
  an~$n$-diagram over a~1-family of cubes.
\end{enumerate}
\end{lm}
\begin{pf} (1) The first part of the statement is immediate. To construct the
  necessary $n$-diagram, group together all the coordinates of $x$
  that are equal and assign the corresponding sets of indices to
  different 0-cubes. This yields an $n$-diagram $\DD$ over a~0-family
  of cubes, and, obviously, $\aut(\DD)$ is exactly the
  $\Ss_n$-stabilizer of $x$ in~$\R^n$.

  (2) Take a nonzero vector $v\,{\in}\,l$. Group together all the equal
  coordinates of $v$, and assign corresponding sets of indices to
  0-cubes, just like we did for $x$. Now, whenever there are two
  groups of coordinates, such that these groups are of equal
  cardinality, and the coordinates in the two groups are negatives of
  each other, we connect the two corresponding 0-cubes with an edge,
  to form a~1-cube. We orient all these cubes in the same coordinate 
  direction.
  Clearly, this yields an $n$-diagram $\DD$ over a~1-family of cubes.
  
  Assume first that our diagram consists of a number of $1$-cubes and
  at most one $0$-cube, with the fiber over this $0$-cube consisting
  of all the indices of the coordinates of~$v$ which are equal to~$0$. 
  The elements of the group $\aut(\DD)$ are of two sorts, depending on
  which of the two elements of $\Z_2$ they are associated to.
  We easily verify that those elements of $\aut(\DD)$, which are
  associated to $0\in\Z_2$, are exactly those $\pi\in\Ss_n$, which
  fix $v$, while those elements of $\aut(\DD)$, which are associated
  to $1\in\Z_2$, are exactly those $\pi\in\Ss_n$, which map $v$
  to~$-v$. Since these are the only two options for mapping $v$, if
  $l$ is to be preserved by the element $\pi$, we have proven the
  lemma in this case.
  
  Assume now that $\DD$ is a diagram of some other form. Then, there
  exist no $\pi\,{\in}\,\Ss_n$ such that $\pi(v)\,{=}\,-v$, i.e., each
  element of $\stab(l)$ fixes $l$ pointwise. In this case,
  $\stab(l)\,{=}\,\stab(v)$, thus we are back to case (1) and the
  diagram can be obtained by splitting all the $1$-cubes into
  $0$-cubes.
\end{pf}


\subsection{Intersections of diagrams}

Let $\CC_1=\{K_1,\dots,K_p\}$, resp.\ $\CC_2=\{L_1,\dots,L_q\}$, be a
$t_1$-, resp.\ $t_2$-family of cubes, where $K_i$ is an $I_1(i)$-cube,
and $L_j$ is an $I_2(j)$-cube, for all $i\in[p]$, $j\in[q]$.

Let $\DD_1$, resp.\ $\DD_2$, be $n$-diagrams over $\CC_1$, resp.\
$\CC_2$:
$$
\begin{array}{ccccl}
\DD_1&:&[n]&\longrightarrow& \vrt(\CC_1),\\
&&k&\mapsto& (K_{\alpha_1(k)},v_k^{(1)}),\\
&&&& \\
\DD_2&:&[n]&\longrightarrow& \vrt(\CC_2),\\
&&k&\mapsto& (L_{\alpha_2(k)},v_k^{(2)}).
\end{array}
$$

\begin{df}
  The {\bf intersection of diagrams} $\DD_1$ and $\DD_2$, denoted
  $\DD=\DD_1\cap\DD_2$, is an $n$-diagram over a $(t_1+t_2)$-family of
  cubes $\CC$ defined as follows:
  $$\CC=\{M_{i,j}\}_{i\in[p],j\in[q]},\quad
  I(i,j)=I_1(i)\cup\{x+t_1\,|\,x\in I_2(j)\},$$
  here $M_{i,j}$ is an $I(i,j)$-cube, furthermore
$$\begin{array}{ccccl}
\DD&:&[n]&\longrightarrow& \vrt(\CC),\\
&&k&\mapsto& (M_{\alpha_1(k),\alpha_2(k)},
v_k^{(1)}\cup\{x+t_1\,|\,x\in v_k^{(2)}\}).
\end{array}$$
\end{df}

Note that the fibers over the vertices and cubes of $\DD$ are
determined by the fibers of $\DD_1$ and $\DD_2$ as follows:
$$\DD^{-1}(M_{i,j})=\DD_1^{-1}(K_i)\cap\DD_2^{-1}(L_j),$$
 and
 \begin{equation} \label{eqfib}
\DD^{-1}(M_{i,j},v)=\DD_1^{-1}(K_i,I_1(i)\cap v)\cap
 \DD_2^{-1}(L_j,\{x-t_1\,|\,x\in v, x>t_1\}),
\end{equation}
for each $v\subseteq I(i,j)$.

\vspace{0.3cm}

\begin{center}
\hspace{-0.6cm}%
  \begin{picture}(0,0)%
    \includegraphics{diag_inter.pstex}%
  \end{picture}%
  \input{diag_inter.pstex_t}%
   \\
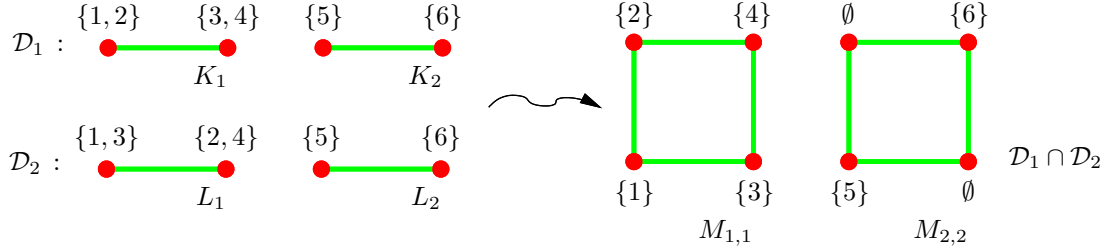
Figure 5. An example of an intersection of two diagrams.
\end{center}

\vspace{0.2cm}

In the above example, observe that $\DD_1\cap\DD_2$ actually contains
two more cubes, $M_{1,2}$ and $M_{2,1}$, with $2$-element index sets
$I(1,2)$ and $I(2,1)$, whose fibers, however, are all empty.

\begin{lm} \label{lm6.6}
  For two $n$-diagrams $\DD_1$ and $\DD_2$, we have
  $\rho(\DD_1\cap\DD_2)=\rho(\DD_1)\wedge\rho(\DD_2),$ where $\wedge$
  denotes the operation of common refinement of the set partitions.
\end{lm}
\begin{pf}
  By \eqref{eqfib}, the blocks of $\rho(\DD_1\cap\DD_2)$ are
  all nonempty intersections of the blocks of $\rho(\DD_1)$ with the
  blocks of $\rho(\DD_2)$, which is precisely the definition of the
  common refinement operation.
\end{pf}

We shall prove two structural theorems about $n$-diagrams. The
first one asserts that taking intersections of diagrams commutes
with passing to the automorphism group.

\begin{thm} \label{thm_dint}
  For two $n$-diagrams $\DD_1$ and $\DD_2$ as above, and
  $\DD=\DD_1\cap\DD_2$ their intersection, we have
 \begin{equation} \label{eq6.3}
\aut(\DD_1)\,\cap\,\aut(\DD_2)\, \, = \, \, \aut(\DD).
 \end{equation}
\end{thm}

\begin{pf}
  First we prove that the set on the left hand side of \eqref{eq6.3}
  is a subset of the set on the right hand side.

  Let $\pi\in\aut(\DD_1)\cap\aut(\DD_2)$. By Definition
  \ref{df6.3} i) we know that $\pi$ preserves the fibers $\DD_1^{-1}(K_i)$,
  for all $i\in[p]$, and $\pi$ preserves the fibers $\DD_2^{-1}(L_j)$, for all
  $j\in[q]$. Hence $\pi$ preserves
  $\DD_1^{-1}(K_i)\cap\DD_2^{-1}(L_j)=\DD^{-1}(M_{i,j})$, for all
  $i\in[p],j\in[q]$, and so property i) of Definition
  \ref{df6.3} is valid for $\pi$.

  By Definition \ref{df6.3} ii), there exist
  $\sigma^{(1)}\in\Z_2^{t_1}$, and $\sigma^{(2)}\in\Z_2^{t_2}$, such
  that $$\sigma_{\alpha_1(k)}^{(1)}(v_k^{(1)})\,= \,v_{\pi(k)}^{(1)},
  \, \, \mbox{ and } \, \,
  \sigma_{\alpha_2(k)}^{(2)}(v_k^{(2)})\,=\, v_{\pi(k)}^{(2)},$$
  for all
  $k\in[n]$, where
  $\sigma_{\alpha_1(k)}^{(1)}=\proj_{I_1(\alpha_1(k))}(\sigma)$, and
  $\sigma_{\alpha_2(k)}^{(2)}=\proj_{I_2(\alpha_2(k))}(\sigma)$.

  Define $\sigma\in\Z_2^{t_1+t_2}$ as a concatenation
  $\sigma=(\sigma^{(1)},\sigma^{(2)})$, that is the first $t_1$
  coordinates of $\sigma$ are equal to $\sigma^{(1)}$, and the last
  $t_2$ coordinates of $\sigma$ are equal to $\sigma^{(2)}$.
  Let $k\in[n]$, and decompose $v_k\subseteq[t_1+t_2]$ as
  $v_k=v_k^{(1)}\cup\tilde v_k^{(2)}$, where
  $v_k^{(1)}=v_k\cap\{1,\dots,t_1\}$, and $\tilde
  v_k^{(2)}=v_k\cap\{t_1+1,\dots,t_1+t_2\}$. Then, we have
\begin{eqnarray*}
\sigma_{\alpha_1(k),\alpha_2(k)}(v_k) & = &
\sigma_{\alpha_1(k),\alpha_2(k)}(v_k^{(1)}\cup\tilde v_k^{(2)}) \\
& = &
\sigma_{\alpha_1(k)}(v_k^{(1)})\cup
\tilde\sigma_{\alpha_2(k)}(\tilde v_k^{(2)})=
v_{\pi(k)}^{(1)}\cup\tilde v_{\pi(k)}^{(2)}=v_{\pi(k)}\, ,
\end{eqnarray*}
where $\sigma_{\alpha_1(k),\alpha_2(k)}=
\proj_{I(\alpha_1(k),\alpha_2(k))}(\sigma)$,
$\tilde\sigma_{\alpha_2(k)}$ is equal to $\sigma_{\alpha_2(k)}$ in the
coordinates $\{t_1+1,\dots,t_1+t_2\}$, and is equal to $0$ in the other
coordinates, while $\tilde v_{\pi(k)}^{(2)}=\{x+t_1\,|\,x\in
v_{\pi(k)}^{(2)}\}$. In other words, $\tilde\sigma_{\alpha_2(k)}$ and
$\tilde v_{\pi(k)}^{(2)}$ are the $t_1$-shifted versions of
$\sigma_{\alpha_2(k)}$ and $v_{\pi(k)}^{(2)}$. So, we have shown that
$\pi\in\aut(\DD_1\cap\DD_2)$.

Now let us prove that the set on the right hand side of \eqref{eq6.3}
is a subset of the set on the left hand side.

Take $\pi\in\aut(\DD)$, then $\pi$ preserves $\DD^{-1}(M_{i,j})$, and
therefore $\pi$ also preserves
\begin{eqnarray*}
\bigcup_{j=1}^q\DD^{-1}(M_{i,j}) & = &
\bigcup_{j=1}^q\DD_1^{-1}(K_i)\cap\DD_2^{-1}(L_j) \, = \,
\DD_1^{-1}(K_i)\cap\bigcup_{j=1}^q\DD_2^{-1}(L_j) \\
& = &
\DD_1^{-1}(K_i)\cap [n]\, = \, \DD_1^{-1}(K_i)\,,
                          \quad \qquad \mbox{ for any }\, i\in[p]\,;
\end{eqnarray*}
in the same way $\pi$
preserves $\DD_2^{-1}(L_j)$, for any $j\,{\in}\,[q]$. This checks
condition i) of Definition~\ref{df6.3}.

Finally, by condition ii) of Definition \ref{df6.3}, there
exists $\sigma\in\Z_2^{t_1+t_2}$, such that for any $k\in[n]$ we have
$\sigma_{\alpha_1(k),\alpha_2(k)}(v_k)=v_{\pi(k)}$. As above, we can
decompose $\sigma=(\sigma^{(1)},\sigma^{(2)})$ and
$v_k=v_k^{(1)}\cup\tilde v_k^{(2)}$ as a~concatenation of the first
$t_1$ and the last $t_2$ coordinates. Then, in the notations which we
used above, we can derive that
$$\sigma_{\alpha_1(k)}^{(1)}(v_k^{(1)})=v_{\pi(k)}^{(1)},\, \, \mbox{ and }
\,\, \, \tilde\sigma_{\alpha_2(k)}^{(2)}(\tilde v_k^{(2)})=\tilde
v_{\pi(k)}^{(2)}.$$
Shifting the second identity down by $t_1$, we get
$\sigma_{\alpha_2(k)}^{(2)}(v_k^{(2)})=v_{\pi(k)}^{(2)}$.
\end{pf}

\subsection{A reduction theorem}

When $\DD$ is an $n$-diagram over a $t$-family of cubes, not every
element $\sigma\in\Z_2^t$ gives rise to an element $\pi\in\aut(\DD)$. The
natural obstruction is that, by Remark~\ref{rem6.4}, fibers with
different cardinalities cannot map to each other. It turns out that
one can always canonically reduce $\DD$ to another $n$-diagram with
the same automorphism group, such that in this new $n$-diagram all 
fibers over vertices in the same cube have the same cardinality.

\begin{thm} \label{thm_red}
  Let $\DD$ be an $n$-diagram over a $t$-family of cubes
  $\CC=(K_1,\dots,K_p)$. Then, there exists an $n$-diagram
  $\widetilde\DD$ over a $\tilde t$-family of cubes
  $\widetilde\CC=(L_1,\dots,L_q)$, such that
\begin{enumerate}
\item[0)] $\ti t\leq t$;
\item[1)] $\aut(\DD)=\aut(\wti\DD)$;
\item[2)] $|\DD^{-1}(L_j,v)|=|\DD^{-1}(L_j,v')|$, for all $j\in [q]$,
  and for all $v,v'\subseteq I_2(j)$, where $I_2(j)$ is the index set
  of~$L_j$.
\end{enumerate}
\end{thm}

In the continuation, we shall call an $n$-diagram satisfying Condition
2) of Theorem~\ref{thm_red} a {\it reduced diagram}.

\vspace{0.3cm}

{\noindent {\bf Proof of Theorem \ref{thm_red}.}}  Let $G$ be the
set of all $\sigma\in\Z_2^t$, such that $\sigma$ occurs as
a~$[t]$-cube symmetry for some $\pi\in\aut(\DD)$. Clearly, $G$ is
a~linear subspace of $\Z_2^t$, when both are viewed as vector spaces
over the field $\F_2$. Hence, there exists $0\leq d\leq t$, such that
$G\cong\Z_2^d$. Therefore, we can choose an~orthogonal linear basis
$\{e_1,\dots,e_t\}$ for $\Z_2^t$, such that $\{e_1,\dots,e_d\}$ is
an~orthogonal linear basis for $G$.

Let us split each cube $K_i\in\CC$ into the orbits of the restriction
of the action of $G$ to $K_i$. We can think of cubes $K_i$ as
coordinate subspaces, that is as intersections of coordinate
hyperplanes, with respect to the standard basis in the vector space
$\Z_2^t$. The orbits themselves however are not coordinate subspaces,
rather they are intersections of the coordinate subspaces
corresponding to cubes with affine linear subspaces of dimension~$d$
obtained from $G$ by parallel translations. Therefore, if we change
the linear basis in $\Z_2^t$ from the standard one to
$\{e_1,\dots,e_t\}$ at the same time as we split the cubes of $\CC$
into the orbits as described above, we end up with a~new $t$-family of
cubes $\wti\CC=(L_1,\dots,L_q)$, and an $n$-diagram $\wti\DD$ over
this family, which is induced from $\DD$.

\vspace{0.3cm}

\begin{center}
  \begin{picture}(0,0)%
    \includegraphics{split.pstex}%
  \end{picture}%
  \input{split.pstex_t}%
  \\[0.2cm]
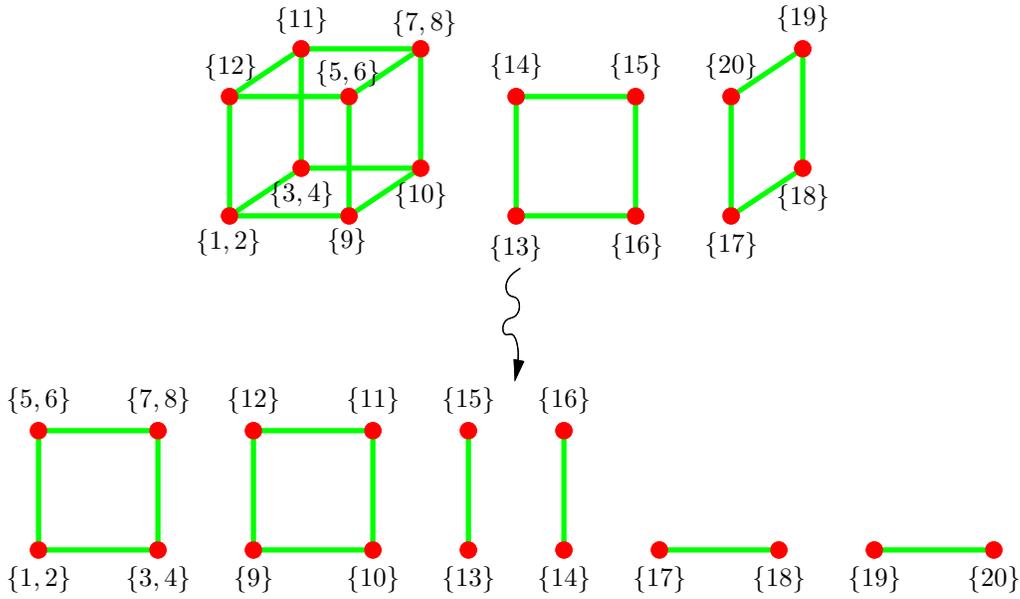
Figure 6. An example of the canonical splitting of a diagram.
\end{center}

\vspace{0.2cm}
By the choice of $G$ and of the basis $\{e_1,\dots,e_t\}$, we see that
all the cubes of $\wti\CC$ actually lie within the coordinate subspace
of $\Z_2^t$ corresponding to the first $d$ coordinates. Thus, we might
as well think of $\wti\CC$ as a~$d$-family of cubes, with $\Z_2^d$
action induced from the action of~$\Z_2^t$, from which
condition~0) of the theorem follows.

Also, since the action on the ground set $[n]$ never changed, we still
have the equality $\aut(\DD)=\aut(\wti\DD)$, verifying
condition~1) of the theorem.

Finally, since $G$ acts transitively on each of its orbits, we can
conclude that the cardinalities of the fibers are constant for the
vertices of the same cube in $\wti\CC$, thus demonstrating the truth
of the last condition, and completing the proof of the theorem.
{\hfill $\Box$\vspace{0.3cm}}


\section{Stabilizers of points in $Y_{\Pi_n}$} \label{sect_stabYP}

In this section we show that the stabilizers of points in
$Y_{\Pi_n}$ are not just abelian, but in fact are isomorphic to direct
products of $\Z_2$. In view of the already proven results, it merely
remains to put the puzzle pieces together.

\begin{thm}
  For $Y_{\Pi_n}$, the De~Concini-Procesi arrangement model of the
  braid arrangement, and $\omega\in Y_{\Pi_n}$, the stabilizer of $\omega$
  with respect to the $\Ss_n$-action on $Y_{\Pi_n}$ is a direct product
  of~$\Z_2$'s:
\[
       \stab_{Y_{\Pi_n}}(\omega)\, \, \cong \, \, \Z_2^h\, , \qquad
  \mbox{ for some }\,  0\leq h\leq\lfloor n/2\rfloor\, .
\]
\end{thm}
\begin{pf}
  By \eqref{eqn_pntYG} a~point in $Y_{\Pi_n}$ can be written as
  $\omega=(x,H_1,l_1,H_2,\dots,H_t,l_t)$, where $H_i\,{\in}\,\Pi_n\,
  {\setminus}\, \{\hat 0\}$,
  and there
  does not exist a subspace $H\,{\in}\,\Pi_n$, $H\,{\neq}\,\R^n$, such that
  $H\,{\supseteq}\,\langle H_t,l_t\rangle$. By Proposition~\ref{prop_stab}
  we know that
\begin{equation}
  \label{eq:6.5.1}
\stab_{Y_{\Pi_n}}(\omega)=\stab_{\R^n}(x)\cap\stab_{\R^n}(l_1)\cap\dots
\cap\stab_{\R^n}(l_t).
\end{equation}
By Lemma~\ref{lm_stab} there exist diagrams $\DD_0,\DD_1,\dots,\DD_t$,
such that
\begin{equation}
  \label{eq:6.5.2}
  \aut(\DD_0)=\stab_{\R^n}(x), \quad \mbox{ and }\,\quad
 \aut(\DD_i)=\stab_{\R^n}(l_i),\, \, \mbox{ for each }\, i\in [t].
\end{equation}

Combining~\eqref{eq:6.5.1}, \eqref{eq:6.5.2}, and
Theorem~\ref{thm_dint}, we find an~$n$-diagram $\DD$, such that
$\aut(\DD)=\stab_{Y_{\Pi_n}}(\omega)$. Moreover, by the Reduction
Theorem~\ref{thm_red}, we can assume that $\DD$ is reduced.

If the partition $\rho(\DD)$ has a~block $B$ of cardinality at
least~3, then, by Lemma~\ref{lm6.6}, so do also the partitions
$\rho(\DD_0),\rho(\DD_1),\dots, \rho(\DD_t)$. Let $H$ be the linear
subspace of $\R^n$ of codimension~2 defined by setting the coordinates
with indices in $B$ equal. By construction, $x\in H$, and
$l_1\subseteq H,\dots,l_t\subseteq H$. Since $H\in\Pi_n$, we see that
$x\in H$ implies $H_1\subseteq H$. Further $l_1\subseteq H$, together
with $H_1\subseteq H$, implies $\langle l_1,H_1\rangle\subseteq H$.
Hence $H_2\subseteq H$, and so on, until we can conclude that $\langle
l_t,H_t\rangle\subseteq H$. This yields a~contradiction, since
$H\neq\R^n$.

So we proved that all blocks of the partition $\rho(\DD)$ are of
cardinality at most~2. Assume now there exist two different blocks
$B_1$ and $B_2$ in $\rho(\DD)$, such that $|B_1|=|B_2|=2$. Let $H$ be
the linear subspace of $\R^n$ of codimension~2 defined by equations
$x_{i_1}=x_{i_2}$, $x_{j_1}=x_{j_2}$, where $B_1=\{i_1,i_2\}$,
$B_2=\{j_1,j_2\}$. Again $H\in\Pi_n$, and by an argument completely
analogous to the previous one, we can trace the two blocks $B_1$ and
$B_2$ through the partitions $\rho(\DD_0),\rho(\DD_1),\dots,
\rho(\DD_t)$, and conclude that $\langle l_t,H_t\rangle\subseteq H$.
This again yields a~contradiction, since $H\neq\R^n$.

Now we know that $\rho(\DD)$ has at most one block of size~2. In
particular, since $\DD$ is reduced, all the fibers over $I$-cubes, for
$|I|\geq 1$, are of cardinality~1. Let us say $\DD$ is an~$n$-diagram
over a~$t$-family of cubes $\CC=(K_1,\dots,K_q)$, where $t$ is minimal
possible. If $\rho(\DD)$ has no blocks of size~2, then there exists
a~group isomorphism between $\aut(\DD)$ and~$\Z_2^t$, since each
element $\pi\in\Z_2^t$ defines the maps between the fibers uniquely.
Each $I$-cube defines at most $|I|$ new directions and has $2^{|I|}$
vertices, hence
$$ t\leq |I_1|+\dots+|I_t|\leq 2^{|I_1|-1}+\dots+2^{|I_t|-1}=n/2.$$

If the partition $\rho(\DD)$ has one block $B$ of size~2, then, since
$\DD$ is reduced, $B$ has to be a~fiber over a~$\emptyset$-cube. With
$t$ chosen as above, it is immediate that
$\aut(\DD)\cong\Z_2^t\times\Z_2$, where the first factor on the right
hand side is the group acting on the $[t]$-cube, and the second factor
is acting on the set~$B$.  Just as before we get $t\leq(n-2)/2$, hence
$t+1\leq n/2$.
\end{pf}



\end{document}